# Les séries congruo-harmoniques alternées
## Partie 2 – Accélérations de convergence

par David Pouvreau[1]

avec les participations de Dominique Tournès[2] et de Gilles Patry[3]

## 1. Introduction

$(p; q) \in (\mathbb{N}^*)^2$ étant fixé, on note $S_{p,q}$ la somme de la série de terme général $\frac{(-1)^k}{pk+q}$, qualifiée de série « congruo-harmonique alternée » (CHA) de paramètres $(p; q)$ :

$$S_{p,q} = \sum_{k=0}^{+\infty} \frac{(-1)^k}{pk+q}.$$

Son explicitation a été déterminée dans Pouvreau (2022). Pour tout $n \in \mathbb{N}$, considérons la somme partielle de rang $n$ de cette série :

$$S_{p,q}^{(n)} = \sum_{k=0}^{n} \frac{(-1)^k}{pk+q}.$$

Un problème majeur est la lenteur de la convergence de la suite $\left(S_{p,q}^{(n)}\right)_{n \in \mathbb{N}}$, qui peut facilement être observée sur les exemples numériques regroupés dans le tableau 1 où, comme dans ceux qui suivront, les décimales correctes sont indiquées en caractères gras :

|  | $n = 100$ | $n = 1\,000$ | $n = 10\,000$ |
| --- | --- | --- | --- |
| $(p; q) = (2; 1)$ | **0,78**7873 | **0,785**647 | **0,7854**2316 |
| $(p; q) = (1; 2)$ | **0,3**11730 | **0,30**7351 | **0,306**90280 |

Tableau 1. Exemples de valeurs des sommes partielles $S_{p,q}^{(n)}$

Il n'est pas superflu de rappeler ici les définitions des notions de « vitesse de convergence » et de « vitesse de convergence relative ». Soit $(u_n)_{n \in \mathbb{N}}$ une suite réelle convergeant vers un réel $\ell$. Si la suite de terme général $\frac{u_{n+1}-\ell}{u_n-\ell}$ est définie et converge vers un réel $V$ de $[0; 1]$, ce nombre $V$ définit la *vitesse de convergence*, ou *taux de convergence*, de $(u_n)_{n \in \mathbb{N}}$.

Si $V = 1$, on dit que $(u_n)_{n \in \mathbb{N}}$ a une convergence lente, ou *infra-linéaire*.

Si $V \in ]0; 1[$, on dit que $(u_n)_{n \in \mathbb{N}}$ converge *linéairement*, ou géométriquement : il existe deux réels $C > 0$ et $k \in ]0; 1[$ tels que $|u_n - \ell| \leq Ck^n$ à partir d'un certain rang. $V$ (appelé taux de convergence linéaire ici) est alors la borne inférieure des nombres $k$ satisfaisant cette condition. Elle peut être déterminée par la preuve de l'existence d'une constante $\lambda$ telle que $|u_n - \ell| \sim_{n \to +\infty} \lambda V^n$. Le nombre de décimales significatives de l'évaluation de $\ell$ par $u_n$ tend ici à croître de manière arithmétique avec $n$, avec une raison qui est une fonction décroissante de $V$.

Enfin, si $V = 0$ on dit que la convergence est *super-linéaire*.

---

[1] Professeur agrégé de mathématiques en CPGE au Lycée Roland Garros, Le Tampon, La Réunion. IREM de La Réunion. Email : pvrdvd@gmail.com.
[2] Professeur de mathématiques à l'Université de La Réunion. Auteur de vérifications numériques, de la figure 3 et de contributions significatives aux démonstrations des théorèmes 2, 3 et 5.
[3] Professeur agrégé de mathématiques et IA-IPR de Rennes. Auteur de la plupart des simulations numériques utilisées pour cet article, effectuées sous Python et vérifiées avec Pari/GP.

Si $(v_n)_{n \in \mathbb{N}}$ est une autre suite convergeant aussi vers $\ell$, on dit qu'elle *converge plus vite* que $(u_n)_{n \in \mathbb{N}}$ si celle de terme général $\left|\frac{v_n - \ell}{u_n - \ell}\right|$ converge vers 0, moins vite si elle diverge vers $+\infty$.

Nous établirons ici que les séries CHA ont effectivement une convergence infra-linéaire. Puis nous construirons plusieurs algorithmes permettant d'accélérer leur convergence. Ces algorithmes sont tous fondés sur le théorème 3 établi dans Pouvreau (2021) :

Pour tout $(p;q) \in (\mathbb{N}^*)^2$ et tout $n \in \mathbb{N}$ :

$$\sum_{k=0}^{+\infty} \frac{(-1)^k}{pk+q} = \sum_{k=0}^{n} \frac{(-1)^k}{pk+q} + \cfrac{(-1)^{n+1}}{\alpha_{p,q}(n) + \cfrac{(1p)^2}{\alpha_{p,q}(n) + \cfrac{(2p)^2}{\alpha_{p,q}(n) + \cfrac{(3p)^2}{\alpha_{p,q}(n) + \cdots}}}} \quad \text{où } \alpha_{p,q}(n) = 2pn + p + 2q.$$

Il s'agira ici d'analyser les effets d'accélération de convergence engendrés par les algorithmes ainsi construits, en précisant notamment la nature de ces convergences.

## 2. Algorithmes ACSCHA d'accélération de la convergence

Pour $n \in \mathbb{N}$ fixé, nous noterons $\rho_{p,q,m}(n)$ la réduite d'ordre $m \in \mathbb{N}$ du développement en fraction continue généralisée de $R_{p,q}(n)$, c'est-à-dire :

$$\rho_{p,q,0}(n) = \frac{1}{\alpha_{p,q}(n)} = \frac{1}{2pn + (p+2q)} \text{ et, pour tout } m \geq 1 : \rho_{p,q,m}(n) = \cfrac{1}{\alpha_{p,q}(n) + \cfrac{(1p)^2}{\alpha_{p,q}(n) + \cfrac{(2p)^2}{\alpha_{p,q}(n) + \cdots + \cfrac{(mp)^2}{\alpha_{p,q}(n)}}}}$$

Le théorème rappelé ci-dessus peut ainsi se reformuler par :

$$\forall (p;q) \in (\mathbb{N}^*)^2, \ \forall n \in \mathbb{N}, \ S_{p,q} = S_{p,q}(n) + (-1)^{n+1} \lim_{m \to +\infty} \rho_{p,q,m}(n).$$

Chaque réduite est un équivalent de $R_{p,q}(n)$ lorsque $n$ est voisin de $+\infty$, dont on peut se servir afin d'accélérer la convergence de la série initiale. Il va s'agir ici de construire avec ces réduites des suites de convergence beaucoup plus rapide que $\left(S_{p,q}(n)\right)_{n \in \mathbb{N}}$.

### 2.1. Vitesse de convergence des séries CHA

Pour tout $(p;q) \in (\mathbb{N}^*)^2$ et tout $n \in \mathbb{N}$ :

$$\left|\frac{S_{p,q}(n+1) - S_{p,q}}{S_{p,q}(n) - S_{p,q}}\right| = \left|\frac{S_{p,q}(n) + \frac{(-1)^{n+1}}{p(n+1)+q} - S_{p,q}}{S_{p,q}(n) - S_{p,q}}\right| = \left|1 - \frac{1}{(pn + (p+q)) R_{p,q}(n)}\right|.$$

Or, il résulte du théorème rappelé en introduction qu'au voisinage de $+\infty$, $R_{p,q}(n) \sim \frac{1}{2pn}$. D'où :

$$\lim_{n \to +\infty} \left|\frac{S_{p,q}(n+1) - S_{p,q}}{S_{p,q}(n) - S_{p,q}}\right| = \left|1 - \frac{1}{p \times \frac{1}{2p}}\right| = 1.$$

### Théorème 1

Toute série CHA converge infra-linéairement vers $S_{p,q}$, avec $\left|S_{p,q}(n) - S_{p,q}\right| \underset{n \to +\infty}{\sim} \frac{1}{2pn}$.

## 2.2. Principe de construction des algorithmes

Observons maintenant, par exemple dans le cas $(p;q) = (2;1)$, les résultats obtenus pour les couples $(n;m)$ correspondant aux premières valeurs de $n$ et de $m$ :

|       | $m=0$   | $m=1$     | $m=2$      | $m=3$       | $m=4$       |
|-------|---------|-----------|------------|-------------|-------------|
| $n=0$ | 0,75000 | 0,8000000 | 0,77777778 | 0,790123457 | 0,782222222 |
| $n=1$ | 0,79167 | 0,7843137 | 0,78571429 | 0,785276074 | 0,785454545 |
| $n=2$ | 0,78333 | 0,7855856 | 0,78536585 | 0,785406302 | 0,785395537 |
| $n=3$ | 0,78631 | 0,7853480 | 0,78540373 | 0,785397206 | 0,785398385 |
| $n=4$ | 0,78492 | 0,7854157 | 0,78539682 | 0,785398328 | 0,785398135 |

Tableau 2 : valeurs de $S_{2,1}{}^{(n)} + (-1)^{n+1}\rho_{2,1,m}{}^{(n)}$ pour $(n;m) \in [\![0;4]\!]^2$

Ce tableau rend manifestes les effets d'accélération engendrés par l'augmentation de $n$ d'une part et par celle de $m$ d'autre part, mais aussi le fait que ces effets se conjuguent si la « direction diagonale » d'une augmentation simultanée de $m$ et $n$ est empruntée.

On peut dès lors envisager d'accélérer d'au moins trois manières la convergence des séries CHA : (1) en se fixant un ordre $m$ de réduction puis en augmentant l'ordre $n$ de sommation partielle ; (2) en se fixant l'ordre $n$ de sommation partielle puis en augmentant l'ordre $m$ de réduction ; (3) en augmentant simultanément les deux ordres $m$ et $n$ (donc avec $m = n$). Nous appellerons ici « ACSCHA » (respectivement U, V et W) les trois algorithmes ainsi construits, et plus généralement tout algorithme construit sur la base du théorème rappelé en introduction.

## 3. Algorithme ACSCHA-U : une accélération de convergence infra-linéaire

Pour tout $(p;q) \in (\mathbb{N}^*)^2$ et tout $m \in \mathbb{N}$, notons $\left(u_{p,q,m}{}^{(n)}\right)_{n\in\mathbb{N}}$ la suite de terme général $S_{p,q}{}^{(n)} + (-1)^{n+1}\rho_{p,q,m}{}^{(n)}$. On l'appellera la suite d'ordre $m$ de l'algorithme ACSCHA-U.

### Théorème 2
<span style="color:red">Toute suite $\left(u_{p,q,m}{}^{(n)}\right)_{n\in\mathbb{N}}$ issue de l'algorithme ACSCHA-U converge infra-linéairement vers $S_{p,q}$ mais plus vite que la suite $\left(S_{p,q}{}^{(n)}\right)_{n\in\mathbb{N}}$, avec :</span>

$$\left|u_{p,q,m}{}^{(n)} - S_{p,q}\right| = \left|\rho_{p,q,m}{}^{(n)} - R_{p,q}{}^{(n)}\right| \sim_{n\to+\infty} \frac{(m+1)!^2}{2^{2m+3}p}\frac{1}{n^{2m+3}}.$$

Pour le montrer, commençons par introduire, pour tout $n \in \mathbb{N}$ provisoirement fixé, les deux suites strictement positives d'entiers naturels $\left(A_{p,q,m}{}^{(n)}\right)_{m\in\mathbb{N}}$ et $\left(B_{p,q,m}{}^{(n)}\right)_{m\in\mathbb{N}}$ telles que :

$$\forall m \in \mathbb{N},\ \begin{cases} A_{p,q,m+2}{}^{(n)} = \alpha_{p,q}{}^{(n)}A_{p,q,m+1}{}^{(n)} + p^2(m+2)^2 A_{p,q,m}{}^{(n)} \\ B_{p,q,m+2}{}^{(n)} = \alpha_{p,q}{}^{(n)}B_{p,q,m+1}{}^{(n)} + p^2(m+2)^2 B_{p,q,m}{}^{(n)} \end{cases}$$

$$\text{avec } \begin{cases} A_{p,q,0}{}^{(n)} = 1 \\ A_{p,q,1}{}^{(n)} = \alpha_{p,q}{}^{(n)} \end{cases} \text{ et } \begin{cases} B_{p,q,0}{}^{(n)} = \alpha_{p,q}{}^{(n)} \\ B_{p,q,1}{}^{(n)} = \left(\alpha_{p,q}{}^{(n)}\right)^2 + p^2 \end{cases}.$$

Ces deux suites satisfont précisément :

$$\forall (m;n) \in \mathbb{N}^2,\ \rho_{p,q,m}{}^{(m)} = \frac{A_{p,q,m}{}^{(n)}}{B_{p,q,m}{}^{(n)}}.$$

Remarquons alors que pour tout $n \in \mathbb{N}$ et tout $m \in \mathbb{N}^*$ :

$$A_{p,q,m+1}{}^{(n)}B_{p,q,m}{}^{(n)} - A_{p,q,m}{}^{(n)}B_{p,q,m+1}{}^{(n)}$$

$$= \left(\alpha_{p,q}{}^{(n)}A_{p,q,m}{}^{(n)} + p^2(m+1)^2 A_{p,q,m-1}{}^{(n)}\right)B_{p,q,m}{}^{(n)} - A_{p,q,m}{}^{(n)}\left(\alpha_{p,q}{}^{(n)}B_{p,q,m}{}^{(n)} + p^2(m+1)^2 B_{p,q,m-1}{}^{(n)}\right)$$

$$= -p^2(m+1)^2\left(A_{p,q,m}{}^{(n)}B_{p,q,m-1}{}^{(n)} - A_{p,q,m-1}{}^{(n)}B_{p,q,m}{}^{(n)}\right).$$

Il en résulte par récurrence finie sur $m$ que :

$$A_{p,q,m+1}{}^{(n)}B_{p,q,m}{}^{(n)} - A_{p,q,m}{}^{(n)}B_{p,q,m+1}{}^{(n)} = (-1)^m p^{2m}(m+1)!^2 \left(A_{p,q,1}{}^{(n)}B_{p,q,0}{}^{(n)} - A_{p,q,0}{}^{(n)}B_{p,q,1}{}^{(n)}\right)$$

$$= (-1)^{m+1} p^{2m+2}(m+1)!^2.$$

Et on peut en déduire :

$$\left|\rho_{p,q,m+1}{}^{(n)} - \rho_{p,q,m}{}^{(n)}\right| = \left|\frac{A_{p,q,m+1}{}^{(n)}}{B_{p,q,m+1}{}^{(n)}} - \frac{A_{p,q,m}{}^{(n)}}{B_{p,q,m}{}^{(n)}}\right| = \left|\frac{A_{p,q,m+1}{}^{(n)}B_{p,q,m}{}^{(n)} - A_{p,q,m}{}^{(n)}B_{p,q,m+1}{}^{(n)}}{B_{p,q,m}{}^{(n)}B_{p,q,m+1}{}^{(n)}}\right|$$

$$= \frac{p^{2m+2}(m+1)!^2}{B_{p,q,m}{}^{(n)}B_{p,q,m+1}{}^{(n)}}.$$

Utilisons alors le fait classique que la suite $\left(\rho_{p,q,m}{}^{(n)}\right)_{m\in\mathbb{N}}$ converge vers $R_{p,q}{}^{(n)}$ avec ses suites extraites de termes de rangs pairs et de rangs impairs qui sont adjacentes. On obtient :

$$\left|\rho_{p,q,m+1}{}^{(n)} - \rho_{p,q,m}{}^{(n)}\right| - \left|\rho_{p,q,m+1}{}^{(n)} - \rho_{p,q,m+2}{}^{(n)}\right| \leq \left|\rho_{p,q,m}{}^{(n)} - R_{p,q}{}^{(n)}\right| \leq \left|\rho_{p,q,m+1}{}^{(n)} - \rho_{p,q,m}{}^{(n)}\right|$$

puis :

$$\frac{p^{2m+2}(m+1)!^2}{B_{p,q,m}{}^{(n)}B_{p,q,m+1}{}^{(n)}} - \frac{p^{2m+4}(m+2)!^2}{B_{p,q,m+1}{}^{(n)}B_{p,q,m+2}{}^{(n)}} \leq \left|\rho_{p,q,m}{}^{(n)} - R_{p,q}{}^{(n)}\right| \leq \frac{p^{2m+2}(m+1)!^2}{B_{p,q,m}{}^{(n)}B_{p,q,m+1}{}^{(n)}}.$$

On a donc aussi :

$$1 - p^2(m+2)^2 \frac{B_{p,q,m}{}^{(n)}}{B_{p,q,m+2}{}^{(n)}} \leq \frac{\left|\rho_{p,q,m}{}^{(n)} - R_{p,q}{}^{(n)}\right|}{\frac{p^2(m+1)!^2}{B_{p,q,m}{}^{(n)}B_{p,q,m+1}{}^{(n)}}} \leq 1.$$

Or, puisqu'il est clair par récurrence sur $m$ que $B_{p,q,m}{}^{(n)} \sim_{n\to+\infty} (2pn)^{m+1}$, on a :

$$\frac{B_{p,q,m}{}^{(n)}}{B_{p,q,m+2}{}^{(n)}} \sim_{n\to+\infty} \frac{1}{(2pn)^2}.$$

D'où le résultat annoncé, l'équivalent obtenu étant évidemment de convergence lente vers 0 :

$$\left|u_{p,q,m}{}^{(n)} - S_{p,q}\right| = \left|\rho_{p,q,m}{}^{(n)} - R_{p,q}{}^{(n)}\right| \sim_{n\to+\infty} \frac{p^{2m+2}(m+1)!^2}{B_{p,q,m}{}^{(n)}B_{p,q,m+1}{}^{(n)}} \sim_{n\to+\infty} \frac{(m+1)!^2}{(2n)^{2m+3}p}.$$

## 4. Algorithme ACSCHA-V : autre accélération de convergence infra-linéaire

On considère ici pour $(p;q) \in (\mathbb{N}^*)^2$ et $n \in \mathbb{N}$ fixés la suite $\left(v_{p,q,m}{}^{(n)}\right)_{m\in\mathbb{N}}$ de terme général $v_{p,q,m}{}^{(n)} = S_{p,q}{}^{(n)} + (-1)^{n+1}\rho_{p,q,m}{}^{(n)}$, appelée suite d'ordre $n$ de l'algorithme ACSCHA-V.

### Théorème 3

Toute suite $\left(v_{p,q,m}{}^{(n)}\right)_{m\in\mathbb{N}}$ issue de l'algorithme ACSCHA-V converge infra-linéairement vers $S_{p,q}$ mais plus vite que la suite $\left(S_{p,q}{}^{(n)}\right)_{n\in\mathbb{N}}$, avec :

$$\exists\, \omega > 0, \ \left|v_{p,q,m}{}^{(n)} - S_{p,q}\right| = \left|\rho_{p,q,m}{}^{(n)} - R_{p,q}{}^{(n)}\right| \sim_{m\to+\infty} \frac{\omega}{m^{2n+1+\frac{2q}{p}}}.$$

Il est d'abord clair que pour tout $n \in \mathbb{N}$ fixé et tout $m \in \mathbb{N}$ :

$$\left|v_{p,q,m}^{(n)} - S_{p,q}\right| = \left|S_{p,q}^{(n)} + (-1)^{n+1}\rho_{p,q,m+1}^{(n)} - S_{p,q}\right| = \left|\rho_{p,q,m}^{(n)} - R_{p,q}^{(n)}\right|.$$

Cherchons des équivalents de $A_{p,q,m}^{(n)}$ et $B_{p,q,m}^{(n)}$. Nous avons vu que ce sont les termes généraux de suites solutions de l'équation de récurrence linéaire à coefficients non constants :

$$\forall m \in \mathbb{N}^*, \quad -C_{p,q,m+1}^{(n)} + \alpha_{p,q}^{(n)} C_{p,q,m}^{(n)} + p^2(m+1)^2 C_{p,q,m-1}^{(n)} = 0.$$

Il s'avère que les deux valeurs propres réelles et toute matrice de passage vers une réduite d'une matrice associée à cette équation sont paramétrées explicitement par $m$. Ce qui entrave l'explicitation d'une suite $\left(C_{p,q,m}^{(n)}\right)_{m \in \mathbb{N}}$. La détermination des équivalents cherchés n'en demeure pas moins possible. En effet, l'équation s'écrit aussi sous la forme :

$$\forall m \in \mathbb{N}^*, \quad -\frac{1}{m^2} C_{p,q,m+1}^{(n)} + \frac{\alpha_{p,q}^{(n)}}{m^2} C_{p,q,m}^{(n)} + \left(p^2 + \frac{2p^2}{m} + \frac{p^2}{m^2}\right) C_{p,q,m-1}^{(n)} = 0.$$

Ce qui montre qu'il s'agit plus précisément d'une équation de récurrence de Poincaré d'ordre 1 qui est irrégulière. La théorie bien établie à cet égard (voir par exemple Hautot (2017), auquel le lecteur est renvoyé ici) permet d'aboutir aux fins visées. L'équation caractéristique associée est :

$$-\frac{1}{m^2} z^2 + \frac{\alpha_{p,q}^{(n)}}{m^2} z + p^2 = 0 \Leftrightarrow -z^2 + \alpha_{p,q}^{(n)} z + p^2 m^2 = 0.$$

Cherchons alors deux solutions asymptotiques sous la forme $z = \rho m^\tau$. On obtient :

$$-\rho^2 m^{2\tau} + \alpha_{p,q}^{(n)} \rho m^\tau + p^2 m^2 = 0.$$

Soit aussi : $p^2 - \rho^2 m^{2\tau-2} + \alpha_{p,q}^{(n)} \rho m^{\tau-2} = 0$. Les seules possibilités (par passage à la limite) se révèlent être $\begin{cases} \tau = 1 \\ \rho = p \end{cases}$ et $\begin{cases} \tau = 1 \\ \rho = -p \end{cases}$. Les évaluations asymptotiques se présentent ainsi sous la forme $\rho m^\tau$, respectivement $r_1 \sim pm^1$ et $r_2 \sim -pm^1$.

Considérons d'abord $r_1$. On pose $C_{p,q,m}^{(n)} = p^m (m-1)! \, \gamma_{p,q,m}^{(n)}$. L'équation de récurrence initiale est alors pour tout $m \geq 2$ équivalente à :

$$-p^{m+1} m! \gamma_{p,q,m+1}^{(n)} + \alpha_{p,q}^{(n)} p^m (m-1)! \gamma_{p,q,m}^{(n)} + p^2 (m^2 + 2m + 1) p^{m-1} (m-2)! \gamma_{p,q,m-1}^{(n)} = 0$$

$$\Leftrightarrow -m(m-1) p^2 \gamma_{p,q,m+1}^{(n)} + \alpha_{p,q}^{(n)} (m-1) p \gamma_{p,q,m}^{(n)} + p^2 (m^2 + 2m + 1) \gamma_{p,q,m-1}^{(n)} = 0$$

$$\Leftrightarrow \left(-1 + \frac{1}{m}\right) \gamma_{p,q,m+1}^{(n)} + \left(\frac{\alpha_{p,q}^{(n)}}{pm} - \frac{\alpha_{p,q}^{(n)}}{pm^2}\right) \gamma_{p,q,m}^{(n)} + \left(1 + \frac{2}{m} + \frac{1}{m^2}\right) \gamma_{p,q,m-1}^{(n)} = 0.$$

Les coefficients de cette équation de récurrence auxiliaire sont :

$$\begin{aligned} a_{2,0} &= -1 \quad ; \quad a_{2,1} = 1 \quad ; \quad a_{2,2} = 0 \\ a_{1,0} &= 0 \quad ; \quad a_{1,1} = \alpha_{p,q}^{(n)}/p \quad ; \quad a_{1,2} = -\alpha_{p,q}^{(n)}/p \\ a_{0,0} &= 1 \quad ; \quad a_{0,1} = 2 \quad ; \quad a_{0,2} = 1. \end{aligned}$$

Ses coefficients caractéristiques sont les nombres déterminés (en convenant que $0^0 = 1$) par :

$$\sigma_{r,s} = \sum_{j=0}^{2} j^r a_{j,s}.$$

Les seuls nécessaires à calculer ici sont $\sigma_{0,0} = 0$ ; $\sigma_{0,1} = 3 + \alpha_{p,q}^{(n)}/p$ ; $\sigma_{1,0} = -2$. La table standard de solutions asymptotiques des équations de Poincaré fournit alors :

$$-\frac{\sigma_{0,1}}{\sigma_{1,0}} = \frac{3 + \frac{\alpha_{p,q}{}^{(n)}}{p}}{2} = \frac{\alpha_{p,q}{}^{(n)} + 3p}{2p} \quad \text{et} \quad \gamma_{p,q,m}{}^{(n)} \underset{m \to +\infty}{\sim} m^{-\frac{\sigma_{0,1}}{\sigma_{1,0}}} = m^{\frac{\alpha_{p,q}{}^{(n)} + 3p}{2p}}.$$

D'où : $C_{p,q,m}{}^{(n)} \underset{m \to +\infty}{\sim} p^m (m-1)!\, m^{\frac{\alpha_{p,q}{}^{(n)} + 3p}{2p}}.$

Cherchons maintenant la solution associée à $r_2$. On pose désormais :

$$C_{p,q,m}{}^{(n)} = (-1)^m p^m (m-1)!\, \gamma_{p,q,m}{}^{(n)}.$$

L'équation de récurrence initiale est ici équivalente pour tout $m \geq 2$ à :

$$\Leftrightarrow \left(-1 + \frac{1}{m}\right) \gamma_{p,q,m+1}{}^{(n)} + \left(-\frac{\alpha_{p,q}{}^{(n)}}{pm} + \frac{\alpha_{p,q}{}^{(n)}}{pm^2}\right) \gamma_{p,q,m}{}^{(n)} + \left(1 + \frac{2}{m} + \frac{1}{m^2}\right) \gamma_{p,q,m-1}{}^{(n)} = 0.$$

On obtient cette fois par la même méthode que précédemment que

$$C_{p,q,m}{}^{(n)} \underset{m \to +\infty}{\sim} (-1)^m p^m (m-1)!\, m^{\frac{3p - \alpha_{p,q}{}^{(n)}}{2p}}.$$

On déduit de cette analyse un système fondamental de solutions asymptotiques de l'équation de récurrence considérée, avec une suite dominée que nous noterons $\left(a_{p,q,m}{}^{(n)}\right)_{m \in \mathbb{N}}$ et une suite dominante que nous noterons $\left(b_{p,q,m}{}^{(n)}\right)_{m \in \mathbb{N}}$, par définition telles que :

$$a_{p,q,m}{}^{(n)} \underset{m \to +\infty}{\sim} (-1)^m p^m (m-1)!\, m^{\frac{3p - \alpha_{p,q}{}^{(n)}}{2p}} \quad ; \quad b_{p,q,m}{}^{(n)} \underset{m \to +\infty}{\sim} p^m (m-1)!\, m^{\frac{3p + \alpha_{p,q}{}^{(n)}}{2p}}.$$

Toute solution de l'équation de récurrence peut s'écrire comme combinaison linéaire de ces deux suites. Par conséquent, en particulier :

$$\exists\, (\lambda_a; \lambda_b; \mu_a; \mu_b) \in \mathbb{R}^4, \ \forall\, m \in \mathbb{N}, \ \begin{cases} A_{p,q,m}{}^{(n)} = \lambda_a a_{p,q,m}{}^{(n)} + \mu_a b_{p,q,m}{}^{(n)} \\ B_{p,q,m}{}^{(n)} = \lambda_b a_{p,q,m}{}^{(n)} + \mu_b b_{p,q,m}{}^{(n)} \end{cases}.$$

Observons enfin que pour $m$ assez grand :

$$\rho_{p,q,m}{}^{(n)} = \frac{A_{p,q,m}{}^{(n)}}{B_{p,q,m}{}^{(n)}} = \frac{\lambda_a a_{p,q,m}{}^{(n)} + \mu_a b_{p,q,m}{}^{(n)}}{\lambda_b a_{p,q,m}{}^{(n)} + \mu_b b_{p,q,m}{}^{(n)}} = \frac{\mu_a + \lambda_a \frac{a_{p,q,m}{}^{(n)}}{b_{p,q,m}{}^{(n)}}}{\mu_b + \lambda_b \frac{a_{p,q,m}{}^{(n)}}{b_{p,q,m}{}^{(n)}}} \underset{m \to +\infty}{\sim} \frac{\mu_a}{\mu_b} = R_{p,q}{}^{(n)}.$$

et qu'on a donc aussi, en posant $\omega = \left|\frac{\lambda_a \mu_b - \lambda_b \mu_a}{\mu_b{}^2}\right|$, qui est non nul du fait de l'indépendance linéaire des suites $\left(A_{p,q,m}{}^{(n)}\right)_{m \in \mathbb{N}}$ et $\left(B_{p,q,m}{}^{(n)}\right)_{m \in \mathbb{N}}$ :

$$\left|\rho_{p,q,m}{}^{(n)} - R_{p,q}{}^{(n)}\right| = \left|\frac{\lambda_a a_{p,q,m}{}^{(n)} + \mu_a b_{p,q,m}{}^{(n)}}{\lambda_b a_{p,q,m}{}^{(n)} + \mu_b b_{p,q,m}{}^{(n)}} - \frac{\mu_a}{\mu_b}\right| = \left|\frac{(\lambda_a \mu_b - \lambda_b \mu_a) a_{p,q,m}{}^{(n)}}{\lambda_b \mu_b a_{p,q,m}{}^{(n)} + \mu_b{}^2 b_{p,q,m}{}^{(n)}}\right|$$

$$\underset{m \to +\infty}{\sim} \left|\frac{\lambda_a \mu_b - \lambda_b \mu_a}{\mu_b{}^2} \cdot \frac{a_{p,q,m}{}^{(n)}}{b_{p,q,m}{}^{(n)}}\right| \underset{m \to +\infty}{\sim} \omega \left|\frac{(-1)^m p^m (m-1)!\, m^{\frac{3p - \alpha_{p,q}{}^{(n)}}{2p}}}{p^m (m-1)!\, m^{\frac{3p + \alpha_{p,q}{}^{(n)}}{2p}}}\right| = \frac{\omega}{m^{\frac{\alpha_{p,q}{}^{(n)}}{p}}}.$$

Disposant désormais d'une évaluation des vitesses de convergence pour les algorithmes U et V, on peut comparer ces vitesses (comparaison dont le sens sera toutefois relativisé plus loin). La suite pertinente à considérer à cet égard est $\left(\left|\frac{u_{p,q,m}{}^{(n)} - S_{p,q}}{v_{p,q,n}{}^{(m)} - S_{p,q}}\right|\right)_{m \in \mathbb{N}}$ :

**Théorème 4**

Si $p > q$, la suite $\left(u_{p,q,m}^{(n)}\right)_{n\in\mathbb{N}}$ converge plus vite que $\left(v_{p,q,n}^{(m)}\right)_{n\in\mathbb{N}}$ pour tout $m \in \mathbb{N}$.

Si $p < q$, la suite $\left(u_{p,q,m}^{(n)}\right)_{n\in\mathbb{N}}$ converge moins vite que $\left(v_{p,q,n}^{(m)}\right)_{n\in\mathbb{N}}$ pour tout $m \in \mathbb{N}$.

Si $p = q$, les suites $\left(u_{p,q,m}^{(n)}\right)_{n\in\mathbb{N}}$ et $\left(v_{p,q,n}^{(m)}\right)_{n\in\mathbb{N}}$ convergent aussi rapidement.

Ce résultat est une conséquence relativement immédiate des théorèmes 2 et 3. En effet :

$$\forall\, m \in \mathbb{N}, \quad \left|\frac{u_{p,q,m}^{(n)} - S_{p,q}}{v_{p,q,n}^{(m)} - S_{p,q}}\right| \sim_{n\to+\infty} \frac{\dfrac{(m+1)!^2}{2^{2m+3}pn^{2m+3}}}{\dfrac{\omega}{n^{2m+1+\frac{2q}{p}}}} = \frac{(m+1)!^2}{2^{2m+3}p\omega} n^{2\left(\frac{q}{p}-1\right)}.$$

D'où la conclusion par le signe de $q/p - 1$, par définition de la vitesse de convergence relative.

## 5. Algorithme ACSCHA-W : une accélération de convergence linéaire

On considère ici la suite $\left(w_{p,q}^{(n)}\right)_{n\in\mathbb{N}}$ de terme général : $w_{p,q}^{(n)} = S_{p,q}^{(n)} + (-1)^{n+1}\rho_{p,q,n}^{(n)}$. On la qualifiera de suite déduite de $\left(S_{p,q}^{(n)}\right)_{m\in\mathbb{N}}$ par l'algorithme ACSCHA-W.

Comme $w_{p,q}^{(n)} = v_{p,q,n}^{(n)} = u_{p,q,n}^{(n)}$ pour tout $n \in \mathbb{N}$, toute suite issue de l'algorithme ACSCHA-W converge plus vite que la suite de mêmes paramètres issue des deux autres algorithmes. Voici quelques exemples de valeurs, spectaculaires eu égard à la petitesse des rangs :

|  | $n = 0$ | $n = 3$ | $n = 5$ | $n = 7$ |
|---|---|---|---|---|
| $w_{1,1}^{(n)}$ | 0,66667 | 0,693146417445 | 0,693147179886527 | 0,693147180559356 |
| $w_{10,1}^{(n)}$ | 0,91667 | 0,938093859970 | 0,938094286672162 | 0,938094287032576 |
| $w_{1,10}^{(n)}$ | 0,05238 | 0,052487740006 | 0,052487740074957 | 0,052487740074975 |

Tableau 3 : quelques valeurs significatives de $w_{p,q}^{(n)}$

En particulier, on obtient vite de la sorte des approximations rationnelles de $\pi$ telles que :

$$\pi \simeq 4w_{2,1}^{(10)} = \frac{3781715948011520}{1203757572990973} \simeq \mathbf{3{,}14159265358979}16$$

Toutes les simulations réalisées suggèrent l'existence d'une propriété universelle :

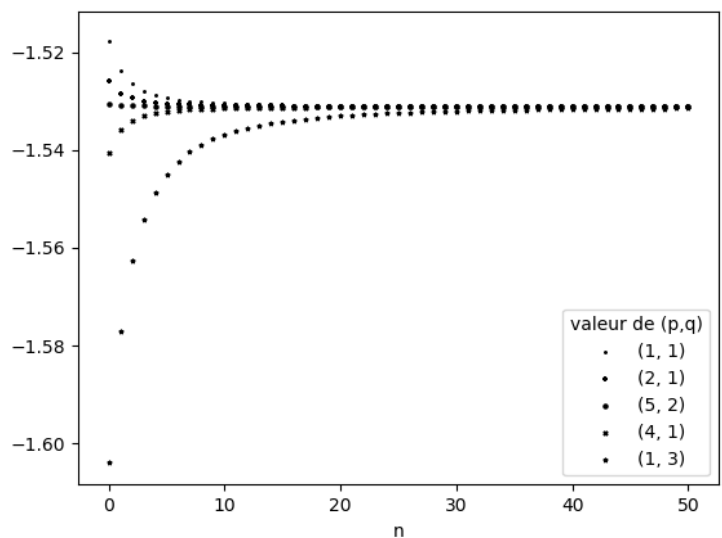

Figure 1 : évolution selon divers choix de $(p; q)$ de $\log\left(\left|\dfrac{w_{p,q}^{(n+1)} - S_{p,q}}{w_{p,q}^{(n)} - S_{p,q}}\right|\right)$ en fonction de $n$

### Conjecture

Toutes les suites issues de l'algorithme ACSCHA-W convergent linéairement à un même taux de convergence linéaire $\chi$ voisin de $0{,}0294372541$, car tel que :

$$\log\left(\left|\frac{w_{4,1}{}^{(1001)} - S_{4,1}}{w_{4,1}{}^{(1000)} - S_{4,1}}\right|\right) \simeq -1{,}531102799 < \log(\chi) < -1{,}531102703 \simeq \log\left(\left|\frac{w_{2,1}{}^{(1001)} - S_{2,1}}{w_{2,1}{}^{(1000)} - S_{2,1}}\right|\right).$$

Si ce résultat n'a pas encore pu être pleinement démontré, on peut au moins établir :

### Théorème 5

Toute suite issue de l'algorithme ACSCHA-W converge linéairement à un taux de convergence linéaire $\chi_{p,q}$ tel que $1/(9e^2) < \chi_{p,q} < 1/(4e^2)$.

La démonstration repose sur un encadrement assez fin des termes de la suite $\left(B_{p,q,m}{}^{(n)}\right)_{m \in \mathbb{N}}$ pour permettre de conclure. On montre d'abord sans difficulté majeure par récurrence double, au moyen de l'identité de récurrence déterminant cette suite, que :

$$\forall\, m \in \mathbb{N}, \qquad \left(\alpha_{p,q}{}^{(n)}\right)^{m+1} \leq B_{p,q,m}{}^{(n)} \leq \left(\alpha_{p,q}{}^{(n)} + pm\right)^{m+1}.$$

Rappelons maintenant l'inégalité établie dans la démonstration du théorème 2 :

$$\forall\, m \in \mathbb{N},\ \forall\, n \in \mathbb{N},\ \left|\rho_{p,q,m}{}^{(n)} - R_{p,q}{}^{(n)}\right| \leq \frac{p^{2m+2}(m+1)!^2}{B_{p,q,m}{}^{(n)} B_{p,q,m+1}{}^{(n)}}.$$

Puisque $B_{p,q,m}{}^{(n)} \geq (2pn)^{m+1}$ pour tout $(m;n) \in \mathbb{N}^2$, on obtient donc, en remarquant par ailleurs que $\left|\rho_{p,q,m}{}^{(n)} - R_{p,q}{}^{(n)}\right| = \left|u_{p,q,m}{}^{(n)} - S_{p,q}\right|$, l'inégalité :

$$\forall\, m \in \mathbb{N},\ \forall\, n \in \mathbb{N},\ \left|u_{p,q,m}{}^{(n)} - S_{p,q}\right| \leq \frac{p^{2m+2}(m+1)!^2}{(2pn)^{m+1}(2pn)^{m+2}} = \frac{(m+1)!^2}{(2n)^{2m+3} p}.$$

On peut en déduire, en l'appliquant à $m = n$ :

$$\forall\, n \in \mathbb{N},\ \left|w_{p,q}{}^{(n)} - S_{p,q}\right| \leq \frac{(n+1)!^2}{(2n)^{2n+3} p}.$$

Il résulte enfin de la formule de Stirling que :

$$\frac{(n+1)!^2}{(2n)^{2n+3} p} \underset{n \to +\infty}{\sim} \frac{\left(\sqrt{2\pi(n+1)}\left(\frac{n+1}{e}\right)^{n+1}\right)^2}{(2n)^{2n+3} p} = \frac{\pi(n+1)\left(\frac{n+1}{e}\right)^{2n+2}}{2^{2n+2} n^{2n+3} p} = \frac{\pi e}{2^{2n+2} e^{2n+3} p}\left(\left(1+\frac{1}{n}\right)^n\right)^2 \left(1+\frac{1}{n}\right)^3$$

$$\underset{n \to +\infty}{\sim} \frac{\pi e}{4^{n+1} e^{2n+3} p} e^2 = \frac{\pi}{4p}\left(\frac{1}{4e^2}\right)^n.$$

Ceci établit que $\left(w_{p,q}{}^{(n)}\right)_{n \in \mathbb{N}}$ converge au moins linéairement et que si sa convergence est linéaire, alors son taux de convergence est inférieur à $\frac{1}{4e^2}$. Par ailleurs, on a aussi pour tout $(m;n) \in \mathbb{N}^2$ :

$$\left|\rho_{p,q,m}{}^{(n)} - \rho_{p,q,m+2}{}^{(n)}\right| \leq \left|\rho_{p,q,m}{}^{(n)} - R_{p,q}{}^{(n)}\right| = \left|u_{p,q,m}{}^{(n)} - S_{p,q}\right|.$$

Or :

$$\left|\rho_{p,q,m}{}^{(n)} - \rho_{p,q,m+2}{}^{(n)}\right| = \left|\frac{A_{p,q,m}{}^{(n)}}{B_{p,q,m}{}^{(n)}} - \frac{A_{p,q,m+2}{}^{(n)}}{B_{p,q,m+2}{}^{(n)}}\right| = \left|\frac{A_{p,q,m}{}^{(n)} B_{p,q,m+2}{}^{(n)} - A_{p,q,m+2}{}^{(n)} B_{p,q,m}{}^{(n)}}{B_{p,q,m}{}^{(n)} B_{p,q,m+2}{}^{(n)}}\right|$$

$$= \left|\frac{\alpha_{p,q}{}^{(n)} p^{2m+2}(m+1)!^2}{B_{p,q,m}{}^{(n)} B_{p,q,m+2}{}^{(n)}}\right| \geq \frac{\alpha_{p,q}{}^{(n)} p^{2m+2}(m+1)!^2}{\left(\alpha_{p,q}{}^{(n)} + pm\right)^{m+1} \left(\alpha_{p,q}{}^{(n)} + p(m+2)\right)^{m+3}}.$$

Pour tout $n \in \mathbb{N}$, on a donc, toujours en utilisant la formule de Stirling :

$$\left|w_{p,q}{}^{(n)} - S_{p,q}\right| \geq \frac{\alpha_{p,q}{}^{(n)} p^{2n+2}(n+1)!^2}{\left(\alpha_{p,q}{}^{(n)} + pn\right)^{n+1} \left(\alpha_{p,q}{}^{(n)} + p(n+2)\right)^{n+3}}$$

$$\sim_{n\to+\infty} \frac{2pn \, p^{2n+2}(n+1)!^2}{(3pn)^{2n+4}} \sim_{n\to+\infty} \frac{(2n)p^{2n+3}}{3^{2n+4}n^{2n+2}n^2 p^{2n+4}} (2\pi)(n+1)\left(\frac{n+1}{e}\right)^{2(n+1)}$$

$$= \frac{4\pi}{9p \, 9^{n+1}e^{2n+2}} \left(\left(1+\frac{1}{n}\right)^n\right)^2 \left(1+\frac{1}{n}\right)^2 \sim_{n\to+\infty} \frac{4\pi}{9p(9e^2)^{n+1}} e^2 = \frac{4\pi}{81p}\left(\frac{1}{9e^2}\right)^n.$$

D'où le fait que $\left(w_{p,q}{}^{(n)}\right)_{n\in\mathbb{N}}$ converge au mieux linéairement à un taux supérieur à $\frac{1}{9e^2}$.

La difficulté cachée dans ce qui reste à prouver de la conjecture tient au fait qu'améliorer l'encadrement de $\chi$ nécessite un encadrement plus fin des termes de $\left(B_{p,q,m}{}^{(n)}\right)_{m\in\mathbb{N}}$ ; et que même ainsi, le problème de l'existence de $\chi$ reste ouvert. Des possibilités d'amélioration sont *a priori* envisageables au moyen d'inégalités faciles à établir par récurrence, telles que :

$$\forall \, m \in \mathbb{N}, \ B_{p,q,m}{}^{(n)} \geq \left(\alpha_{p,q}{}^{(n)}\right)^{m+1} + \frac{m(m+1)(2m+1)}{6} p^2 \left(\alpha_{p,q}{}^{(n)}\right)^{m-1}.$$

Mais l'inégalité précédente conduit par exemple à une majoration de $\left|w_{p,q}{}^{(n)} - S_{p,q}\right|$ par une suite d'équivalent $\frac{36\pi}{pn^2}\left(\frac{1}{4e^2}\right)^n$ ; ce qui ne permet pas d'améliorer l'évaluation du taux de convergence.

La nature de la convergence des suites issues des algorithmes U, V et W est illustrée sur la figure 2 du point de vue du nombre de décimales dites « significatives », ou « précision », de l'approximation de $S_{p,q}$ qu'on obtient avec leurs termes respectifs. C'est-à-dire au moyen de la partie entière des termes des suites $\left(-\log\left(\left|u_{p,q,m}{}^{(n)} - S_{p,q}\right|\right)\right)_{n\in\mathbb{N}}$, $\left(-\log\left(\left|v_{p,q,n}{}^{(m)} - S_{p,q}\right|\right)\right)_{n\in\mathbb{N}}$ et $\left(-\log\left(\left|w_{p,q}{}^{(n)} - S_{p,q}\right|\right)\right)_{n\in\mathbb{N}}$. La différence entre les accélérations induites par U et V prévue par le théorème 4 en cas de nette dissymétrie entre $p$ et $q$ peut aussi être discernée dans cet exemple :

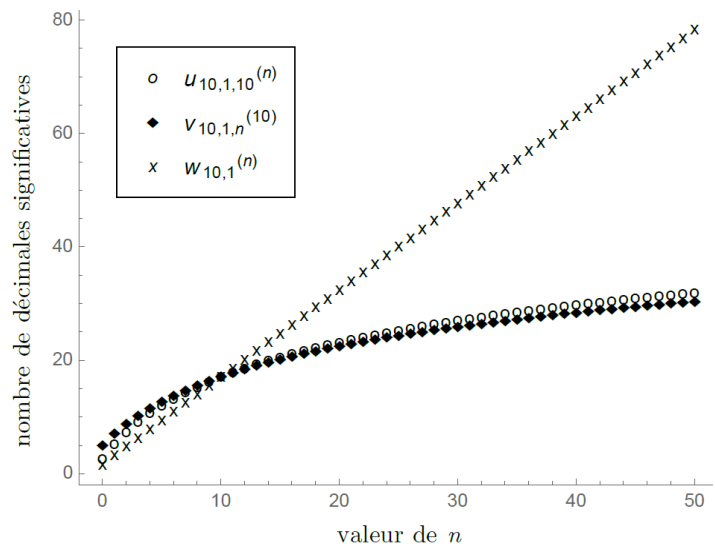

Figure 2

## 6. Algorithmes ACSCHA d'accélération de convergence super-linéaire

Existe-t-il des algorithmes ACSCHA engendrant des convergences super-linéaires ? On peut songer à en construire par extractions de la suite $\left(w_{p,q}^{(n)}\right)_{n\in\mathbb{N}}$ : sa très rapide convergence linéaire s'y prête naturellement. Néanmoins, ce n'est pas nécessaire. Une sorte de « semi-extraction » de $\left(w_{p,q}^{(n)}\right)_{n\in\mathbb{N}}$ au moyen d'une application « semi-extractrice » $\zeta$ opérant sur le seul ordre de sommation partielle, qui implique moins d'opérations (avec certes une moindre précision) que si $\zeta$ opère simultanément sur l'ordre de réduction, se révèle en effet suffisante pour obtenir une convergence super-linéaire dès lors que $\zeta$ est de croissance assez rapide :

### Théorème 6
Soit $\zeta : \mathbb{N} \to \mathbb{N}$ une application strictement croissante. Si $n = o_{n\to+\infty}(\zeta(n))$, alors la suite $\left(w_{\zeta,p,q}^{(n)}\right)_{n\in\mathbb{N}}$ de terme général $S_{p,q}^{(\zeta(n))} + (-1)^{\zeta(n)+1}\rho_{p,q,n}^{(\zeta(n))}$, dite issue de l'algorithme ACSCHA-$W_\zeta$, converge super-linéairement, avec :

$$\left|w_{\zeta,p,q}^{(n)} - S_{p,q}\right| \sim_{n\to+\infty} \frac{\pi}{4p}\left(\frac{1}{4e^2}\right)^n \left(\frac{n}{\zeta(n)}\right)^{2n+3}.$$

La preuve de ces résultats repose sur les inégalités établies en section 5 :

$$\forall\, m \in \mathbb{N},\ \forall\, n \in \mathbb{N},\ \frac{\alpha_{p,q}^{(n)} p^{2m+2}(m+1)!^2}{\left(\alpha_{p,q}^{(n)} + pm\right)^{m+1} \left(\alpha_{p,q}^{(n)} + p(m+2)\right)^{m+3}} \leq \left|u_{p,q,m}^{(n)} - S_{p,q}\right| \leq \frac{(m+1)!^2}{(2n)^{2m+3}p}$$

En résulte :

$$\forall\, n \in \mathbb{N},\ \frac{\alpha_{p,q}^{(\zeta(n))} p^{2n+2}(n+1)!^2}{\left(\alpha_{p,q}^{(\zeta(n))} + pn\right)^{n+1} \left(\alpha_{p,q}^{(\zeta(n))} + p(n+2)\right)^{n+3}} \leq \left|w_{\zeta,p,q}^{(n)} - S_{p,q}\right| \leq \frac{(n+1)!^2}{(2\zeta(n))^{2n+3}p}.$$

En utilisant là encore la formule de Stirling, on obtient de part et d'autre :

$$\frac{(n+1)!^2}{(2\zeta(n))^{2n+3}p} \sim_{n\to+\infty} \frac{2\pi(n+1)\left(\frac{n+1}{e}\right)^{2n+2}}{2^{2n+3}p(\zeta(n))^{2n+3}} \sim_{n\to+\infty} \frac{\pi}{4pe^2} n^{2n+3}\left(\left(1+\frac{1}{n}\right)^n\right)^2 \frac{1}{(2e)^{2n}} \frac{1}{(\zeta(n))^{2n+3}}$$

$$\sim_{n\to+\infty} \frac{\pi}{4p}\left(\frac{1}{4e^2}\right)^n \left(\frac{n}{\zeta(n)}\right)^{2n+3}.$$

$$\frac{\alpha_{p,q}^{(\zeta(n))} p^{2n+2}(n+1)!^2}{\left(\alpha_{p,q}^{(\zeta(n))} + pn\right)^{n+1} \left(\alpha_{p,q}^{(\zeta(n))} + p(n+2)\right)^{n+3}} \sim_{n\to+\infty} \frac{2p\zeta(n)p^{2n+2} 2\pi n\left(\frac{n+1}{e}\right)^{2n+2}}{(2p\zeta(n))^{2n+4}}$$

$$\sim_{n\to+\infty} \frac{\pi}{4pe^2} \frac{n^{2n+3}}{(4e^2)^n} \frac{1}{(\zeta(n))^{2n+3}}\left(\left(1+\frac{1}{n}\right)^n\right)^2 \sim_{n\to+\infty} \frac{\pi}{4p}\left(\frac{1}{4e^2}\right)^n \left(\frac{n}{\zeta(n)}\right)^{2n+3}.$$

D'où effectivement : $\left|w_{\zeta,p,q}^{(n)} - S_{p,q}\right| \sim_{n\to+\infty} \frac{\pi}{4p}\left(\frac{1}{4e^2}\right)^n \left(\frac{n}{\zeta(n)}\right)^{2n+3}$.

La nature super-linéaire de la convergence de $\left(w_{\zeta,p,q}^{(n)}\right)_{n\in\mathbb{N}}$ si $n = o_{n\to+\infty}(\zeta(n))$ en découle par :

$$\left|\frac{w_{\zeta,p,q}^{(n+1)} - S_{p,q}}{w_{\zeta,p,q}^{(n)} - S_{p,q}}\right| \sim_{n\to+\infty} \frac{1}{4e^2}\left(\frac{n+1}{\zeta(n+1)}\right)^{2n+5}\left(\frac{\zeta(n)}{n}\right)^{2n+3} \sim_{n\to+\infty} \frac{1}{4}\left(\frac{n}{\zeta(n+1)}\right)^2\left(\frac{\zeta(n)}{\zeta(n+1)}\right)^{2n+3} < \frac{1}{4}\left(\frac{n}{\zeta(n)}\right)^2.$$

Une illustration numérique de ce résultat : avec $\zeta : n \to n^2$, le nombre $4w_{\zeta,2,1}{}^{(10)}$ donne une approximation de $\pi$ à la précision 37, tandis que $4w_{2,1}{}^{(10)}$ le fait seulement à la précision 16. La figure 3 illustre la comparaison des cas $\zeta : n \to n^3$ et $\zeta : n \to 2^n$ avec la suite $(w_{2,1}{}^{(n)})$.

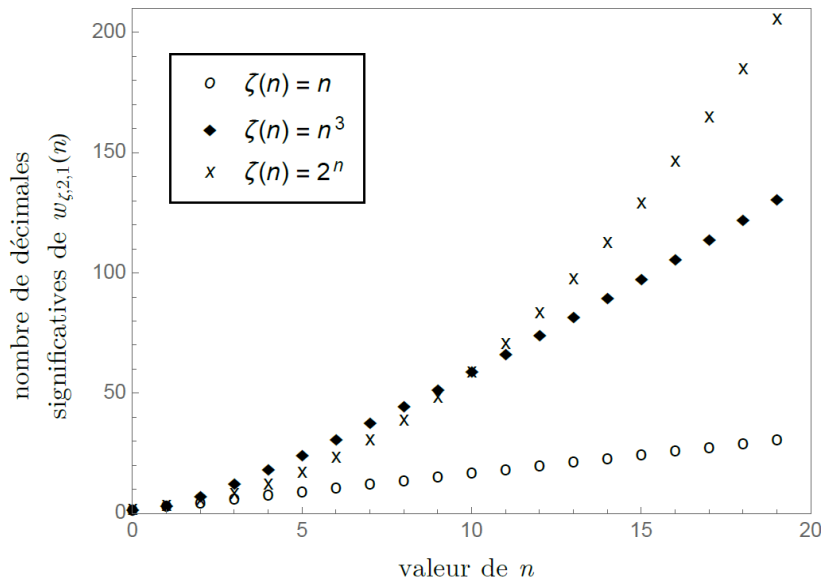

Figure 3

La question des conditions pour qu'une suite $(w_{p,q}{}^{(\zeta(n))})_{n\in\mathbb{N}}$ extraite de $(w_{p,q}{}^{(n)})_{n\in\mathbb{N}}$ converge super-linéairement demeure toutefois pertinente. En effet, si le théorème 6 montre (avec les précédents) que tout choix de $\zeta$ dominant l'identité permet d'engendrer ce type de convergence pour $(w_{p,q}{}^{(\zeta(n))})_{n\in\mathbb{N}}$, il reste à préciser la condition la plus faible sur $\zeta$ permettant ce résultat. Une illustration de la rapidité de convergence de ces suites extraites : avec $\zeta : n \to n^2$, $4w_{2,1}{}^{(\zeta(5))}$ suffit pour réaliser comme $4w_{\zeta,2,1}{}^{(10)}$ l'approximation de $\pi$ à la précision 37.

## 7. Problèmes d'optimalité et de « performance » relative des algorithmes

Nous avons comparé les *vitesses de convergence* de suites issues d'algorithmes ACSCHA. Néanmoins, une comparaison pertinente des *algorithmes* impose de prendre en compte non seulement ces vitesses, mais aussi le nombre d'opérations élémentaires nécessaires pour obtenir des résultats voisins. Il n'existe toutefois pas, dans la littérature scientifique, de critère standard concernant le « coût » opératoire, son calcul dépendant de conventions et du point de vue choisi (mathématique, informatique, énergétique, etc.). Pas plus qu'il n'existe de critère commun pour évaluer la « performance » relative de deux algorithmes d'accélération de convergence intégrant adéquatement ces deux aspects de la vitesse et du « coût » opératoire.

Plusieurs problèmes demeurent donc ouverts à cet égard. Comment définir un concept précis de « performance » relative qui soit opérationnel en particulier pour l'étude comparée des algorithmes U, V, W et $W_\zeta$ ? Il n'est pas *a priori* certain, du point de vue esquissé ici, que W soit plus « performant » que U et V, et que $W_\zeta$ le soit plus que W ; et ce, même en dépit des différences de nature entre leurs convergences respectives. On peut aussi envisager d'examiner de ce point de vue un problème lié à celui évoqué plus haut concernant la comparaison entre les suites issues de $W_\zeta$ et les suites extraites de celles issues de W : existe-t-il pour toute extractrice $\zeta$ donnée dominant l'identité une autre extractrice $\theta$ dominée par $\zeta$ telle que la suite extraite $(w_{p,q}{}^{(\theta(n))})_{n\in\mathbb{N}}$ converge super-linéairement tout en étant au moins aussi « performante » que $(w_{\zeta,p,q}{}^{(n)})_{n\in\mathbb{N}}$ ?

Les problèmes qui viennent d'être évoqués relèvent en fait aussi d'une problématique intimement liée : celle de l'optimalité algorithmique. Existe-t-il un algorithme ACSCHA ou une classe de tels algorithmes qui soi(en)t au moins aussi « performant(s) » que tous les autres ? On peut notamment s'interroger sur l'éventualité que W soit une solution de ce problème, quitte à inclure dans la discussion les suites extraites de suites issues de cet algorithme. L'importance à cet égard du rôle de W est notamment suggérée par les considérations suivantes, certes indépendantes des questions de « coût » calculatoire. Fixons $N \in \mathbb{N}^*$ pair, puis faisons varier $n$ dans $[\![1;N]\!]$. Des simulations comme celle illustrée par la figure 4 suggèrent que le nombre de décimales exactes de $S_{p,q}$ obtenues avec $S_{p,q}^{(n)} + (-1)^{n+1}\rho_{p,q,m}^{(n)}$ est maximal pour $n = N/2$.

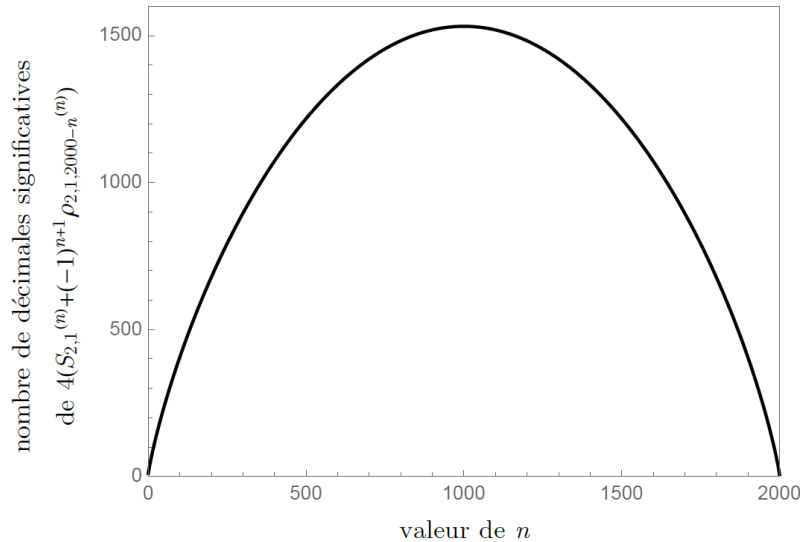

Figure 4

Une diversification des simulations laisse en fait plutôt penser que l'optimalité pour $m = n$ n'est qu'asymptotiquement vraie, c'est-à-dire pour $N$ voisin de $+\infty$. Plus précisément, le maximum apparaît toujours pour $n$ voisin de $N/2$ si $p > q$, mais ceci semble n'être vrai lorsque $q > p$ que si $N$ est grand. Une forme d'optimalité de W apparaît donc au moins en ce sens, qui reste à établir.

## 8. Comparaisons avec d'autres algorithmes d'accélération de convergence

On peut mettre en perspective les algorithmes ACSCHA avec d'autres plus classiques tels que les transformations d'Euler et de Van Wijngaarden, les évaluations asymptotiques de la fonction de Hurwitz-Lerch, l'algorithme de Newton-Cotes ou les approximations sommatoires de Boole. Au moins dans ces deux derniers cas, les algorithmes ACSCHA se révèlent plus efficaces en termes d'accélération de convergence. C'est aussi le cas concernant le delta 2 d'Aitken :

### Théorème 7

$\left(u_{p,q,0}^{(n)}\right)_{n\geq 1}$ est la suite déduite de $\left(S_{p,q}^{(n)}\right)_{n\geq 2}$ par l'algorithme du delta-2 d'Aitken. Celui-ci engendre donc pour les séries CHA une accélération dans l'ensemble significativement moindre que les algorithmes ACSCHA.

En effet, dans le cas de la série CHA de paramètres $(p;q)$, on obtient pour tout $n \geq 2$ :

$$S_{p,q}^{(n-1)} = S_{p,q}^{(n)} - \frac{(-1)^n}{pn+q}$$

$$S_{p,q}{}^{(n-2)} = S_{p,q}{}^{(n)} - (-1)^n \left( \frac{1}{pn+q} - \frac{1}{pn+(q-p)} \right) = S_{p,q}{}^{(n)} + (-1)^n \frac{p}{(pn+q)(pn+(q-p))}.$$

La suite déduite de $\left(S_{p,q}{}^{(n)}\right)$ par l'algorithme du delta-2 d'Aitken a donc pour terme général :

$$\forall n \geq 2, \ A_{p,q}{}^{(n)} = \frac{S_{p,q}{}^{(n)} S_{p,q}{}^{(n-2)} - \left(S_{p,q}{}^{(n-1)}\right)^2}{\left(S_{p,q}{}^{(n)} - S_{p,q}{}^{(n-1)}\right) - \left(S_{p,q}{}^{(n-1)} - S_{p,q}{}^{(n-2)}\right)}$$

$$= \frac{S_{p,q}{}^{(n)} \left( S_{p,q}{}^{(n)} + (-1)^n \frac{p}{(pn+q)(pn+(q-p))} \right) - \left( S_{p,q}{}^{(n)} - \frac{(-1)^n}{pn+q} \right)^2}{\frac{(-1)^n}{pn+q} - \frac{(-1)^{n-1}}{pn+(q-p)}}$$

$$= \frac{\frac{S_{p,q}{}^{(n)}}{pn+q}\left(2 + \frac{p}{pn+(q-p)}\right) - \frac{(-1)^n}{(pn+q)^2}}{\frac{1}{pn+q} + \frac{1}{pn+(q-p)}} = S_{p,q}{}^{(n)} + (-1)^{n+1} \frac{pn+(q-p)}{(pn+q)(2pn+(2q-p))}$$

$$= S_{p,q}{}^{(n-1)} + (-1)^{n+1} \left( \frac{pn+(q-p)}{(pn+q)(2pn+(2q-p))} - \frac{1}{pn+q} \right)$$

$$= S_{p,q}{}^{(n-1)} + \frac{(-1)^n}{2pn+(2q-p)} = S_{p,q}{}^{(n-1)} + (-1)^n \rho_{p,q,0}{}^{(n-1)} = u_{p,q,0}{}^{(n-1)}.$$

## 9. Conclusion

L'évaluation exacte des restes partiels des séries CHA par des fractions continues généralisées a permis d'élaborer plusieurs algorithmes d'accélération de convergence ACSCHA. L'étude de ces algorithmes s'est révélée riche du point de vue de l'analyse ; elle a nécessité la connexion de multiples outils, dont la théorie des équations de récurrence de Poincaré.

Il est établi que les suites issues des algorithmes U et V convergent infra-linéairement et d'autant plus rapidement que l'ordre de ces suites est élevé. La dépendance directe au signe de $p - q$ de la vitesse de convergence relative de deux suites respectivement issues de ces deux algorithmes a été mise en évidence. Il est aussi acquis que celles issues de l'algorithme W ont une très rapide convergence linéaire de taux compris entre $1/(9e^2)$ et $1/(4e^2)$. Des simulations numériques poussées ont suscité la conjecture digne de confiance que ce taux est une constante $\chi$ universelle relativement à $(p;q)$, dont une valeur approximative de précision 10 a été donnée. Ce résultat reste à démontrer. Il s'est enfin révélé possible de construire une classe infinie d'algorithmes ACSCHA (dits $W_\zeta$) dont les suites ont une convergence super-linéaire, ce qui est remarquable au vu de l'extrême lenteur de la convergence des séries initialement considérées. Un problème associé qui demeure ouvert est celui des conditions les plus faibles possibles pour qu'une suite extraite d'une suite issue de W converge super-linéairement.

La nécessité est toutefois apparue de construire un instrument général de mesure de la « performance » relative de deux algorithmes d'accélération de convergence, qui prenne en compte à la fois la vitesse de convergence et son « coût » opératoire.  Seul un tel instrument permettrait en effet une comparaison vraiment pertinente des algorithmes ACSCHA, rendant alors non moins pertinente la question de savoir s'il en existe de « performance » maximale – d'autant plus que l'algorithme W présente au moins certains traits d'optimalité.

Toute cette étude a été nourrie par un dialogue constructif entre simulation numérique et construction théorique, analogue à celui qui s'observe entre physique théorique et physique expérimentale. La simulation n'a jamais été possible que sur des bases et questionnements théoriques préalables. Elle a en retour suscité de nouvelles questions dont certaines demeurent ouvertes et permis de conjecturer l'existence de propriétés *a priori* peu évidentes que la théorie s'est ensuite efforcée d'établir, en fournissant parfois même à celle-ci des indications utiles pour l'élaboration des démonstrations ; et pour, au final, en confirmer les résultats expérimentalement.

**Bibliographie**